\documentclass[a4paper,11pt]{article}
\usepackage{amsmath}
\usepackage{amsfonts}
\usepackage{amssymb}
\usepackage{amscd}
\usepackage{pb-diagram}
\usepackage{epic}

\newtheorem{prop}{Proposition}

\newtheorem{teo}{Theorem}
\newtheorem{corol}{Corollary}
\date{}
\begin{document}

%%%%%%%%%%%%%%%%%%%%%%%%%%%%%%%%%%%%%%%%%%%%%%%%%%%%%%%%%%%%%%%%%%%

\title{Corwin-Greenleaf multiplicity function for compact\\
extensions of the Heisenberg group %\footnote{Typeset title in 10~pt Times Roman uppercase and boldface. Please write down in pencil a short title to be used as the running head.}
}

\author{Majdi Ben Halima and Anis Messaoud}

\maketitle

\begin{abstract} Let $\mathbb{H}_n$ be the
$(2n+1)$-dimensional Heisenberg group and $K$ a closed subgroup of
$U(n)$ acting on $\mathbb{H}_n$ by automorphisms such that
$(K,\mathbb{H}_n)$ is a Gelfand pair. Let $G=K\ltimes\mathbb{H}_n$
be the semidirect product of $K$ and $\mathbb{H}_n$. Let
$\mathfrak{g}\supset\mathfrak{k}$ be the respective Lie algebras
of $G$ and $K$, and $\operatorname{pr}:
\mathfrak{g}^{*}\longrightarrow\mathfrak{k}^{*}$ the natural
projection. For coadjoint orbits
$\mathcal{O}^{G}\subset\mathfrak{g}^{*}$ and
$\mathcal{O}^{K}\subset\mathfrak{k}^{*}$, we denote by
$n\big(\mathcal{O}^{G},\mathcal{O}^{K}\big)$ the number of
$K$-orbits in $\mathcal{O}^{G}\cap \operatorname{pr}^{-1}(\mathcal{O}^{K})$,
which is called the Corwin-Greenleaf multiplicity function. In
this paper, we give two sufficient conditions on $\mathcal{O}^G$
in order that
$$n\big(\mathcal{O}^G,\mathcal{O}^K\big)\leq 1\:\:\text{for any $K$-coadjoint
orbit}\:\:\mathcal{O}^{K}\subset\mathfrak{k}^{*}.$$ For $K=U(n)$,
assuming furthermore that $\mathcal{O}^{G}$ and $\mathcal{O}^{K}$
are admissible and denoting respectively by $\pi$ and $\tau$ their
corresponding irreducible unitary representations, we also discuss
the relationship between $n\big(\mathcal{O}^G,\mathcal{O}^K\big)$
and the multiplicity $m(\pi,\tau)$ of $\tau$ in the restriction of
$\pi$ to $K$. Especially, we study in Theorem 4 the case where $n(\mathcal{O}^{G},\mathcal{O}^{K})\neq m(\pi,\tau)$. This inequality is interesting because we expect the equality as the naming of the Corwin-Greenleaf multiplicity function suggests.
\end{abstract}

{\footnotesize{\textbf{Keywords.}Heisenberg motion group, generic unitary representation, generic coadjoint orbit, Corwin-–Greenleaf multiplicity function.}
\vspace{0,2cm}\\
\textbf{Mathematics Subject Classification 2000.}
22E20, 22E45, 22E27, 53C30}

%%%%%%%%%%%%%%%%%%%%%%%%%%%%%%%%%%%%%%%%%%%%%%%%%%%%%%%%%%%%%%%%%%%%%%%%%%%%%%%%%%%%%%%%%%%%%
\section{Introduction}
%%%%%%%%%%%%%%%%%%%%%%%%%%%%%%%%%%%%%%%%%%%%%%%%%%%%%%%%%%%%%%%%%%%%%%%%%%%%%%%%%%%%%%%%%%%%%
Let $G$ be a connected and simply connected nilpotent Lie group
with Lie algebra $\mathfrak{g}$ and $\widehat{G}$ the unitary dual of $G$, i.e. the set of all
equivalence classes of irreducible unitary representations of $G$. Then Kirillov proved 
that the unitary dual $\widehat{G}$ of $G$ is parametrized by
$\mathfrak{g}^{*}/G$, the set of coadjoint orbits. The bijection
$$\widehat{G}\simeq\mathfrak{g}^{*}/G$$ is called the Kirillov
correspondence (see [7]). Let $\pi$ be the unitary representation corresponding to a given coadjoint orbit $\mathcal{O}^{G}\subset\mathfrak{g}^{*}$. Let $K$ be a subgroup of $G$. Then the restriction $\pi\big\vert_{K}$ is decomposed into a direct integral of irreducible unitary representations of $K$:
$$\pi\big\vert_{K}\simeq \displaystyle \int_{\widehat{K}}^{\oplus} m(\pi,\tau)d\mu(\tau) \hspace{0.5cm} \text{(branching rule)}$$
where $d\mu$ is a Borel measure on the unitary dual $\widehat{K}$. Then Corwin and Greenleaf proved that the above multiplicity $m(\pi,\tau)$ coincides almost everywhere with the ``mod $K$'' intersection number $n(\mathcal{O}^{G},\mathcal{O}^{K})$ defined as follows:
$$n(\mathcal{O}^{G},\mathcal{O}^{K}):=\sharp\Big(\big(\mathcal{O}^{G}\cap \operatorname{pr}^{-1}(\mathcal{O}^{K})\big)/K\Big)$$
(see [4]). Here, $\mathcal{O}^{G}\subset\mathfrak{g}^{*}$ and $\mathcal{O}^{K}\subset\mathfrak{k}^{*}$ are the coadjoint orbits corresponding to $\pi\in\widehat{G}$ and $\tau\in\widehat{K}$, respectively, under the Kirillov correspondence $\widehat{G}\simeq\mathfrak{g}^{*}/G$ and $\widehat{K}\simeq\mathfrak{k}^{*}/K$, and $$\operatorname{pr}:\mathfrak{g}^{*}\longrightarrow\mathfrak{k}^{*}$$ is the natural projection. The function
$$n:\mathfrak{g}^{*}/G\times\mathfrak{k}^{*}/K\longrightarrow\mathbb{N}\cup\{\infty\},\:\:\:\:\:(\mathcal{O}^{G},\mathcal{O}^{K})\longmapsto n(\mathcal{O}^{G},\mathcal{O}^{K})$$ is sometimes referred as the \textbf{Corwin-Greenleaf multipliccity function}. In the special case that $\tau=1_K$, the formula for the multiplicity function $n(\mathcal{O}^{G},\left\lbrace 0\right\rbrace )$ is
$$n(\mathcal{O}^{G},\left\lbrace 0\right\rbrace ):=\sharp\Big(\big(\mathcal{O}^{G}\cap \mathfrak{k}^{\perp}\big)/K\Big),$$ where $\mathfrak{k}^{\perp}:=\operatorname{pr}^{-1}(\left\lbrace 0\right\rbrace )=\left\lbrace \ell\in\mathfrak{g}^*; \ell(\mathfrak{k})=0\right\rbrace   $. \vspace{0,2cm}

In the spirit of the orbit method due to Kirillov, R. Lipsman established a bijection between a class of coadjoint orbits of $G$ and the unitary dual $\widehat{G}$ (see [13]). Given a linear form $\psi\in\mathfrak{g}^{*}$, we denote by $G(\psi)$ its stabilizer in $G$. Then $\psi$ is called admissible if there exists a unitary character $\chi$ of the identity component of $G(\psi)$ such that $d\chi=i\psi_{\vert\mathfrak{g}(\psi)}$. Let $\mathfrak{g}^{\ddagger}$ be the set of all admissible linear forms on $\mathfrak{g}$. For $\psi\in\mathfrak{g}^{\ddagger}$, one can construct an irreducible unitary repesentation $\pi_{\psi}$ by holomorphic induction. According to Lipsman [13], every irreducible unitary representation of $G$ arises in this manner. By observing that $\pi_{\psi}$ is equivalent to $\pi_{\psi^{\,'}}$ if and only if $\psi$ and $\psi^{\,'}$ lie in the same $G$-orbit, we get finally a bijection between the space $\mathfrak{g}^{\ddagger}/G$ of admissible coadjoint orbits and $\widehat{G}$.\vspace{0,2cm}

Let $\pi \in\widehat{G}$ and $\tau\in\widehat{K}$ correspond to admissible coadjoint orbits $\mathcal{O}^{G}$ and $\mathcal{O}^{K}$ respectively and let $\operatorname{pr}:\mathfrak{g}^{*}\longrightarrow\mathfrak{k}^{*}$ be the restriction map. One expects that the multiplicity of $\tau$ in $\pi\big\vert_{K}$ is given by $\sharp\Big(\big(\mathcal{O}^{G}\cap \operatorname{pr}^{-1}(\mathcal{O}^{K})\big)/K\Big)$. Results in this direction have been established for compact extensions of $\mathbb{R}^n$ (see [1]). In this setting the Corwin-Greenleaf multiplicity function $n(\mathcal{O}^{G},\mathcal{O}^{K})$ may become greater than one, or even worse, may take infinity. For example, if $(K,\mathbb{H}_n)$ is a Gelfand pair then $n(\mathcal{O}^{G},\left\lbrace 0\right\rbrace ) = 1$, i.e., $\mathcal{O}^{G}\cap\mathfrak{k}^{\perp}$ is a single $K$-orbit (see [3]).\vspace{0,2cm}

\textbf{Question.} Give a sufficient condition on the admissible coadjoint orbit $\mathcal{O}^G$ in $\mathfrak{g}^*$ in order that
$$n(\mathcal{O}^{G},\mathcal{O}^{K}) \leq 1\:\:\text{for any admissible coadjoint orbit}\:\: \mathcal{O}^K\subset\mathfrak{k}^*.$$
Our interest for this question is motivated by the formulation and the results by Kobayashi-Nasrin [9,16] which may be intepreted as the ``classical limit'' of the multiplicity-free theorems in the branching laws of semisimple Lie groups that were established in [10,11,12] by three different methods, explicit branching laws [10], the theory of visible actions [11], and Verma modules [12].\vspace{0,2cm}

Let $\mathbb{H}_n=\mathbb{C}^n\times\mathbb{R}, n\geq 1$, be the standard Heisenberg group of real dimension $2n+1$. The maximal compact subgroup of $Aut(\mathbb{H}_n)$ is the unitary group $U(n)$, and it acts by $k.(z,t)=(kz,t)$. In this paper we consider the Lie group $G=K\ltimes\mathbb{H}_n$, the semidirect product of the $K$ and $\mathbb{H}_n$, where $K$ stands for a closed subgroup of $U(n)$ acting on $\mathbb{H}_n$ as above. Our group $G$ is obviously a subgroup of the so-called Heisenberg motion group, which is the semidirect product $U(n)\ltimes\mathbb{H}_n$. 
The group $K$ acts on the unitary dual $\widehat{\mathbb{H}_n}$ of $\mathbb{H}_n$ via $$ k.\sigma=\sigma\circ k^{-1}$$ for $k\in K$ and $\sigma\in\widehat{\mathbb{H}_n}$. Let $K_{\sigma}$ denote the stabilizer of $\sigma$ (up to unitary equivalence). Let $\pi$ be an irreducible unitary representation of $G$ associated to a given admissible coadjoint orbit $\mathcal{O}$ in $\mathfrak{g}^{\ddagger}/G$. Mackey's little group theory [14,15] tells us that $\pi$ is determined by a pair $(\sigma,\tau)$ where $\sigma\in\widehat{\mathbb{H}_n}$ and $\tau\in\widehat{K_{\sigma}}$. We consider here the case where the representation $\pi$ is generic, i.e., $\pi$ has Mackey parameters $(\sigma,\tau)$ such that the stabilizer $K_{\sigma}$ is all of $K$. In this case we have $$\pi(k,z,t)=\tau(k)\otimes\sigma(z,t)\circ W_{\sigma}(k),$$ $(k,z,t)\in G$, with $W_{\sigma}$ being a (non-projective) unitary representation of $K$ in the Hilbert space $\mathcal{H}_{\sigma}$ of $\sigma$ that intertwines $k.\sigma$ with $\sigma$: $$(k.\sigma)(z,t)=W_{\sigma}(k)^{-1}\circ\sigma(z,t)\circ W_{\sigma}(k)$$ for all $k\in K, (z,t)\in\mathbb{H}_n$. The main results of
the present work are\vspace{0,2cm}\\
\textbf{Theorem 1.} If $(K,\mathbb{H}_n)$ is a Gelfand pair and $U$ is a
central element of $\mathfrak{k}$, then
$$n\big(\mathcal{O}_{(U,0,x)}^G,\mathcal{O}^{K}_{X}\big)\leq 1$$
for any coadjoint orbit $\mathcal{O}_{X}^{K}$ in
$\mathfrak{k}^{*}$.\vspace{0,2cm}\\
\textbf{Theorem 2.} We have $$m(\pi_{(\lambda,\alpha)},\tau_{\mu})\neq 0\Rightarrow n(\mathcal{O}^{G}_{(\lambda,\alpha)},\mathcal{O}^{K}_{\mu})\neq
0.$$
\textbf{Theorem 3.} Let $n\geq 2$. Assume that $\lambda$ is strongly dominant weight of $K=U(n)$. Then for any dominant weight $\mu$ of $K$ such that $B_{\lambda,\mu}$ is invertible we have $$n(\mathcal{O}^{G}_{(\lambda,\alpha)},\mathcal{O}^{K}_{\mu})\leq 1.$$The matrix $B_{\lambda,\mu}$ is defined in Section 4.3 p 11.\vspace{0,2cm}\\
\textbf{Theorem 4.} Let $n\geq 2$. If the dominant weight $\lambda=(\lambda_1,...,\lambda_n)$ of $K$ satisfies $\lambda_1=...=\lambda_n=a$ for some $a\in\mathbb{Z}$, then for any dominant weight $\mu$ of $K$ with $\mu\neq\lambda$ we have
$$n(\mathcal{O}^{G}_{(\lambda,\alpha)},\mathcal{O}^{K}_{\mu})\leq 1$$ Moreover, $n(\mathcal{O}^{G}_{(\lambda,\alpha)},\mathcal{O}^{K}_{\mu})\neq0$ if and only if $\mu$ is of the form\vspace{0.2cm}\\
\textbf{Case 1:} if $\alpha>0$ then $\mu=(\underbrace{b,...,b}_p,\underbrace{a,...,a}_q)\in\mathbb{Z}^n$, $p+q =n$, $b\in \mathbb{Z}$ with $b>a.$\\
\textbf{Case 2:} if $\alpha<0$ then $\mu=(\underbrace{a,...,a}_p,\underbrace{b,...,b}_q)\in\mathbb{Z}^n$, $p+q =n$, $b\in \mathbb{Z}$ with $a>b.$\\
Consequently, if $\mu_{n-1}\neq a$ and $n(\mathcal{O}^{G}_{(\lambda,\alpha)},\mathcal{O}^{K}_{\mu})\neq 0$ then $m(\pi_{(\lambda,\alpha)},\tau_{\mu})\neq n(\mathcal{O}^{G}_{(\lambda,\alpha)},\mathcal{O}^{K}_{\mu})$.
\vspace{0,2cm}

The paper is organized as follows. Section 2 introduces the coadjoint orbits of $K\ltimes \mathbb{H}_n$. In Sec. 3, we give two sufficient conditions on $\mathcal{O}^G$ in order that $n\big(\mathcal{O}^G,\mathcal{O}^K\big)\leq 1$ for any $K$-coadjoint orbit $\mathcal{O}^{K}\subset\mathfrak{k}^{*}$. Section 4.1 deals with the description of the generic unitary dual $\widehat{U(n)\ltimes \mathbb{H}_n}$ of $U(n)\ltimes \mathbb{H}_n$. Section 4.2 is devoted to the description of the subspace of generic admissible coadjoint orbits of $U(n)\ltimes \mathbb{H}_n$ and to the branching rules from $U(n)\ltimes \mathbb{H}_n$ to $U(n)$. In Sec. 4.3, the Corwin-Greenleaf multiplicity function for $U(n)\ltimes \mathbb{H}_n$ is studied in some situations and the main results of this work are derived.

%%%%%%%%%%%%%%%%%%%%%%%%%%%%%%%%%%%%%%%%%%%%%%%%%%%%%%%%%%%%%%%%%%%%%%%%%%%%%%%%%%%%%%%%%%%
\section{Coadjoint orbits of $K\ltimes\mathbb{H}_{n}$}
%%%%%%%%%%%%%%%%%%%%%%%%%%%%%%%%%%%%%%%%%%%%%%%%%%%%%%%%%%%%%%%%%%%%%%%%%%%%%%%%%%%%%%%%%%%%%%%
On the $n$-dimensional complex vector space $\mathbb{C}^n$, we fix
the usual scalar product $\langle \:,\:\rangle$. Let
$\mathbb{H}_n=\mathbb{C}^n\times\mathbb{R}$ with group law
$$(z,t)(z',t'):=(z+z',t+t'-\frac{1}{2}Im\langle z,z'\rangle)$$
denote the $(2n+1)$-dimensional Heisenberg group. Let $K$ be a
closed subgroup of $U(n)$. The group $K$ acts naturally on
$\mathbb{H}_{n}$ by automorphisms, and then one can form the
semidirect product $G=K\ltimes\mathbb{H}_{n}$. Let us denote by
$(k,z,t)$ the elements of $G$ where $k\in K$ and
$(z,t)\in\mathbb{H}_n$. The group law of $G$ is given by

$$(k,z,t)\cdot(k',z',t')=(kk',z+kz',t+t'-\frac{1}{2}Im\langle z,kz'\rangle).$$
We identify the Lie algebra $\mathfrak{h}_n$ of $\mathbb{H}_n$
with $\mathbb{H}_n$ via the exponential map. We also identify the
Lie algebra $\mathfrak{k}$ of $K$ with its vector dual space
$\mathfrak{k}^*$ through the $K$-invariant inner product
$$(A,B)= tr(AB).$$
For $z\in\mathbb{C}^n$ define the $\mathbb{R}$-linear form $z^*$
in $(\mathbb{C}^n)^*$ by $$z^*(w):=Im\langle z,w\rangle.$$ One
defines a map $\times:
\mathbb{C}^n\times\mathbb{C}^n\longrightarrow\mathfrak{k}$,
$(z,w)\mapsto z\times w$ by $$(z\times w,B)=z\times w(B):=w^*(Bz)$$ with
$B\in\mathfrak{k}$. It is easy to verify that for $k\in K$, one
has
$$\operatorname{Ad}_{{}_{K}}(k)(z\times w)=(kz)\times(kw).$$
Each element $\nu$ in
$\mathfrak{g}^*=(\mathfrak{k}\ltimes\mathfrak{h}_n)^*$ can be
identified with an element $(U,u,x)\in
\mathfrak{k}\times\mathbb{C}^n\times\mathbb{R}$ such that
$$\langle (U,u,x), (B,w,s)\rangle=(U, B) + u^*(w) + xs,$$ where $(B,w,s)\in\mathfrak{g}$. By a direct computation, one obtains that the coadjoint action of $G$ is
$$\operatorname{Ad}_{{}_{G}}^*(k,z,t)(U,u,x)=\big(\operatorname{Ad}_{{}_K}(k)U+z\times
(ku)+\frac{x}{2}z\times z,ku+xz,x\big).$$ Letting $k$ and $z$ vary
over $K$ and $\mathbb{C}^n$ respectively, the coadjoint orbit
$\mathcal{O}_{(U,u,x)}^G$ through the linear form $(U,u,x)$ can be
written
$$\mathcal{O}_{(U,u,x)}^G=\Big\{\big(\operatorname{Ad}_{{}_K}(k)U+z\times
(ku)+\frac{x}{2}z\times z,ku+xz,x\big);\: k\in K,
z\in\mathbb{C}^n\Big\}$$ or equivalently, replacing $z$ by $kz$,
$$\mathcal{O}_{(U,u,x)}^G=\Big\{k\cdot\big(U+z\times
u+\frac{x}{2}z\times z,u+xz,x\big);\: k\in K,
z\in\mathbb{C}^n\Big\}.$$ \textbf{Remark} Here we regard $z$ as a column vector $z=(z_1,...z_n)^T$ and $z^*:=\overline{z}^T$. Then $z\times u\in\mathfrak{u}^*(n)\cong\mathfrak{u}(n)$ is the $n$ by $n$ skew Hermitian matrix $\frac{i}{2}(uz^*+zu^*)$. Indeed, for all $B\in\mathfrak{u}(n)$ we compute$$(uz^*+zu^*,B)=tr((uz^*+zu^*)B)=\sum_{1\leq i,j\leq n}B_{ji}z_i\overline{u}_j - \sum_{1\leq i,j\leq n}u_i\overline{B}_{ij}\overline{z}_j = -2iz\times u(B).$$In particular, $z\times z$ is the skew Hermitian matrix $izz^*$ whose entries are determined by $(izz^*)_{lj}=iz_l\overline{z}_j$.\\

The $G$-coadjoint orbit arising from the
initial point $(U,0,x)(x\neq 0)$ is said to be generic. Notice that the
space of generic coadjoint orbits of $G$ is parametrized by the
set
$\big(\mathfrak{k}/K\big)\times\big(\mathbb{R}\setminus\{0\}\big)$.
Concluding this section, let us underline that the union of all
generic coadjoint orbits of $G$ is dense in $\mathfrak{g}^*$.
%%%%%%%%%%%%%%%%%%%%%%%%%%%%%%%%%%%%%%%%%%%%%%%%%%%%%%%%%%%%%%%%%%%%%%%%%%%%%%%%%%%%%%%%%%%%%%%%%%%%%%%%%%%%%%%%%
\section{Corwin-Greenleaf multiplicity function for $K\ltimes\mathbb{H}_n$}
%%%%%%%%%%%%%%%%%%%%%%%%%%%%%%%%%%%%%%%%%%%%%%%%%%%%%%%%%%%%%%%%%%%%%%%%%%%%%%%%%%%%%%%%%%%%%%%%%%%%%%%%%%%%%%%%%%%%%%%%%%%%%%
We keep the notation of Sec. 2. Consider the generic coadjoint
orbit $\mathcal{O}_{(U,0,x)}^G$ through the element $(U,0,x)$ in
$\mathfrak{g}^{*}$. For $X\in\mathfrak{k}$, we introduce the set
$$\mathcal{F}_{X}:=\Big\{z\in\mathbb{C}^n;\:U+\frac{x}{2}z\times
z\in\mathcal{O}_{X}^{K}\Big\}.$$ Here $\mathcal{O}^{K}_{X}$ is the
$K$-coadjoint orbit in $\mathfrak{k}^{*}\simeq\mathfrak{k}$
through $X$. Letting $H$ be the stabilizer of $U$ in $K$, we
define an equivalence relation in $\mathcal{F}_X$ by
$$z\sim w\Leftrightarrow \exists h\in H;\:w=hz.$$ The set of equivalence classes is denoted by
$\mathcal{F}_{X}/H$.
\begin{prop} For any $X\in\mathfrak{k}$, we have
$$n\big(\mathcal{O}_{(U,0,x)}^G,\mathcal{O}^{K}_{X}\big)=\sharp\Big(\mathcal{F}_{X}/H\Big).$$
\end{prop}
\textbf{Proof.} Fix an element $X$ in $\mathfrak{k}$. For
$z\in\mathbb{C}^n$, let us set
$$E_z:=\Big\{k\cdot\big(U+\frac{x}{2}z\times
z,xz,x\big);\: k\in K\Big\}.$$ Observe that
$$E_{z}=E_{w}\Leftrightarrow z\sim w.$$ Since
$$\mathcal{O}_{(U,0,x)}^G\cap
\operatorname{pr}^{-1}(\mathcal{O}^{K}_{X})=\underset{z\in\mathcal{F}_X}{\bigcup}E_{z},$$
it follows that
\begin{eqnarray*}
n\big(\mathcal{O}_{(U,0,x)}^G,\mathcal{O}^{K}_{X}\big)&=&\sharp\Big[\big(\mathcal{O}_{(U,0,x)}^G\cap
\operatorname{pr}^{-1}(\mathcal{O}^{K}_{X})\big)/K\Big]\\&=&\sharp\Big(\mathcal{F}_{X}/H\Big).
\end{eqnarray*}
This completes the proof of the proposition.$\hfill\square$\\

Following [2], we define the moment map
$\tau:\mathbb{C}^n\rightarrow\mathfrak{k}^{*}$ for the natural
action of $K$ on $\mathbb{C}^n$ by
$$\tau(z)(A)=z^{*}(Az)$$
for $A\in\mathfrak{k}$. Since $\langle z,Az\rangle$ is pure
imaginary, one can also write $\tau(z)(A)=\frac{1}{i}\langle
z,Az\rangle$. The map $\tau$ is a key ingredient in the proof of
the following result.
\begin{teo} If $(K,\mathbb{H}_n)$ is a Gelfand pair and $U$ is a
central element of $\mathfrak{k}$, then
$$n\big(\mathcal{O}_{(U,0,x)}^G,\mathcal{O}^{K}_{X}\big)\leq 1$$
for any coadjoint orbit $\mathcal{O}_{X}^{K}$ in
$\mathfrak{k}^{*}$.
\end{teo}
\textbf{Proof.} Let $U$ be a central element of $\mathfrak{k}$.
Then for any $X\in\mathfrak{k}$,
$$n\big(\mathcal{O}_{(U,0,x)}^G,\mathcal{O}^{K}_{X}\big)=\sharp\Big(\mathcal{F}_{X}/K\Big).$$
Fix a non-zero element $X\in\mathfrak{k}$ and assume that the set
$\mathcal{F}_X$ is not empty. It is clear that $\mathcal{F}_X$ is
stable under the natural action of $K$ on $\mathbb{C}^n$. If $z$
and $w$ are two elements in $\mathcal{F}_X$, then there exists
$k\in K$ such that
$$w\times w=\operatorname{Ad}_{{}_{K}}(k)(z\times z).$$
Thus we get the equality
$\mathcal{O}_{\tau(z)}^{K}=\mathcal{O}_{\tau(w)}^{K}$. Since
$(K,\mathbb{H}_n)$ is a Gelfand pair, the moment map
$\tau:\mathbb{C}^{n}\rightarrow\mathfrak{k}^{*}$ is injective on
$K$-orbits [2]. That is, if
$\mathcal{O}_{\tau(z)}^{K}=\mathcal{O}_{\tau(w)}^{K}$, then
$Kz=Kw$. We conclude that the $K$-action on $\mathcal{F}_X$ is
transitive and hence
$n\big(\mathcal{O}_{(U,0,x)}^G,\mathcal{O}^{K}_{X}\big)=1$.$\hfill\square$\\

%%%%%%%%%%%%%%%%%%%%%%%%%%%%%%%%%%%%%%%%%%%%%%%%%%%%%%%%%%%%%%%%%%%%%%%%%%%%%%%%%%%%%%%%%%%%%%%%%%%%%%%%%%%%%%%%%%%%%%%%%%%%%%%%%%%%%%%%%%%%%%%%%%%%%%%%%%%%%%%%%%%%%%%
\section{Corwin-Greenleaf multiplicity function for
$U(n)\ltimes\mathbb{H}_n$ and branching rules}
%%%%%%%%%%%%%%%%%%%%%%%%%%%%%%%%%%%%%%%%%%%%%%%%%%%%%%%%%%%%%%%%%%%%%%%%%%%%%%%%%%%%%%%%%%%%%%%%%%%%%%%%%%%%%%%%%%%%%%%%%%%%%%%%%%%%%%%%%%%%%%%%%%%%%%%%%%%%%%%%%%%%%%
%We shall freely use the notation of the previous sections. Let us
%fix $K=U(n)$ with $n\geq 1$. Then $G=K\ltimes\mathbb{H}_n$ is
%so-called Heisenberg motion group... \vspace{1cm}

%Next, we shall explicitly determine all the admissible coadjoint
%orbits $\mathcal{O}^{K}\in\mathfrak{k}^{*}$ satisfying
%$n(\mathcal{O}^G,\mathcal{O}^K)\neq 0$.

%%%%%%%%%%%%%%%%%%%%%%%%%%%%%%%%%%%%%%%%%%%%%%%%%%%%%%%%%%%%%%%%%%%%%%%%%%%%%%%%%%%%%%%%%%%
\subsection{Generic unitary dual of $U(n)\ltimes\mathbb{H}_{n}$}
%%%%%%%%%%%%%%%%%%%%%%%%%%%%%%%%%%%%%%%%%%%%%%%%%%%%%%%%%%%%%%%%%%%%%%%%%%%%%%%%%%%%%%%%%%%%%%%
In the sequel, we fix $K=U(n)$ with $n\geq 2$. Then
$G=K\ltimes\mathbb{H}_n$ is the so-called Heisenberg motion group. The description of the unitary dual $\widehat{G}$ of $G$ is based on the Mackey little group theory. In the present paper we consider only the generic irreducible unitary representation of $G$.\vspace{0.2cm}

Let us recall a useful fact from the representation theory of the
Heisenberg group (see,e.g.,[5] for details). The infinite dimensional irreducible
representations of $\mathbb{H}_n$ are parametrized by
$\mathbb{R}^*$. For each $\alpha\in\mathbb{R}^*$, the Kirillov
orbit $\mathcal{O}_\alpha^{\mathbb{H}_n}$ of the irreducible
representation $\sigma_\alpha$ is the hyperplane
$\mathcal{O}_\alpha^{\mathbb{H}_n}=\{(z,\alpha),\:z\in\mathbb{C}^n\}$.
It is clear that for every $\alpha$ the coadjoint orbit
$\mathcal{O}_\alpha$ is invariant under the $K$-action. Therefore
$K$ preserves the equivalence class of $\sigma_\alpha$. The
representation $\sigma_\alpha$ can be realized in the Fock space
$$\mathcal{F}_\alpha(n)=\Big\{f:\mathbb{C}^n\longrightarrow\mathbb{C}\:\text{holomorphic}\:| \int_{\mathbb{C}^n}\vert f(w)\vert^2 e^{-\frac{\vert\alpha\vert}{2}\vert
w\vert^2} dw<\infty\Big\}$$ as
\begin{equation*}
\sigma_\alpha(z,t)f(w)=e^{i\alpha
t-\frac{\alpha}{4}|z|^2-\frac{\alpha}{2}\langle
 w,z\rangle}f(w+z)
\end{equation*}
for $\alpha>0$ and
\begin{equation*}
\sigma_\alpha(z,t) f(\overline{w})=e^{i\alpha
t+\frac{\alpha}{4}|z|^2+\frac{\alpha}{2}\langle
 \overline{w},\overline{z}\rangle}f(\overline{w}+\overline{z})
\end{equation*}
for $\alpha<0$. We refer the reader to [5] or [6] for a discussion of the Fock space. For each $A\in K$, the operator
$W_\alpha(A):\mathcal{F}_{\alpha}(n)\rightarrow
\mathcal{F}_{\alpha}(n)$ defined by
$$
W_\alpha(A)f(w)=f(A^{-1}w)
$$
intertwines $\sigma_\alpha$ and $(\sigma_{\alpha})_A$ given by $(\sigma_{\alpha})_A(z,t):=\sigma_{\alpha}(Az,t)$. Observe that
$W_\alpha$ is a unitary representation of $K$ in the Fock space
$\mathcal{F}_{\alpha}(n)$.\vspace{0,2cm}

As usual, the dominant weights of $K=U(n)$ are parametrized by
sequences $\lambda=(\lambda_1,...,\lambda_n)\in\mathbb{Z}^n$ such
that $\lambda_1\geq\lambda_2\geq...\geq\lambda_{n}$. Denote by
$(\tau_{\lambda},\mathcal{H}_{\lambda})$ an irreducible unitary
representation of $K$ with highest weight $\lambda$. Then by Mackey [15], for each
nonzero $\alpha\in\mathbb{R}$
\begin{equation*}
\pi_{(\lambda,\alpha)}(A,z,t):=\tau_\lambda(A)\otimes\sigma_\alpha(z,t)\circ
W_\alpha(A),\hspace{0.2cm} (A,z,t)\in G,
\end{equation*}
is an irreducible unitary representation of $G$ realized in
$\mathcal{H}_{\lambda}\otimes\mathcal{F}_{\alpha}(n)$. This
representation $\pi_{(\lambda,\alpha)}$ is said to be generic. The
set of all equivalence classes of generic irreducible unitary
representations of $G$, denoted by $\widehat{G}_{gen}$, is called
the generic unitary dual of $G$. Notice that $\widehat{G}_{gen}$
has full Plancherel measure in the unitary dual $\widehat{G}$ (see [8]).
%%%%%%%%%%%%%%%%%%%%%%%%%%%%%%%%%%%%%%%%%%%%%%%%%%%%%%%%%%%%%%%%%%%%%%%%%%%%%%%%%%%%%%%%%%%%%%%%%%%%%%%%%%%%%%%%%%%%%%%%%%%
\subsection{Generic admissible coadjoint orbits of $U(n)\ltimes\mathbb{H}_{n}$ and Branching rules}
%%%%%%%%%%%%%%%%%%%%%%%%%%%%%%%%%%%%%%%%%%%%%%%%%%%%%%%%%%%%%%%%%%%%%%%%%%%%%%%%%%%%%%%%%%%%%%%%%%%%%%%%%%%%%%%%%%%%%%%%%% 
We shall freely use the notation of the previous subsection. Given
a dominant weight $\lambda=(\lambda_{1},...,\lambda_{n})$ of $K$,
we associate to $\pi_{(\lambda,\alpha)}$ the linear form
$\ell_{\lambda,\alpha}=(U_\lambda,0,\alpha)$ in $\mathfrak{g}^*$
where
$$U_\lambda=
\left( \begin{array}{ccc}
i\lambda_1& \ldots & 0\\
\vdots & \ddots & \vdots\\
0 & \ldots & i\lambda_n
\end{array}\right).$$
Observe that $\ell_{\lambda,\alpha}$ is an admissible linear form
on $\mathfrak{g}$. Denote by $G\big(\ell_{\lambda,\alpha}\big)$,
$K\big(\ell_{\lambda,\alpha}\big)$ and
$\mathbb{H}_n\big(\ell_{\lambda,\alpha}\big)$ the stabilizers of
$\ell_{\lambda,\alpha}$ respectively in $G$, $K$ and
$\mathbb{H}_n$.  We have
\begin{eqnarray*}
G\big(\ell_{\lambda,\alpha}\big)&=& \{(A,z,t)\in G;\:(AU_{\lambda}A^* +\frac{\alpha}{2}z\times z,\alpha z,\alpha)=(U_\lambda,0,\alpha)\} \\
&=& \{(A,0,t)\in G; AU_{\lambda}A^*=U_{\lambda}\},\\
K\big(\ell_{\lambda,\alpha}\big)&=& \{A\in K;\:(AU_{\lambda}A^*,0,\alpha)=(U_{\lambda},0,\alpha)\} \\
&=& \{A\in K;\:AU_{\lambda}A^*=U_{\lambda}\},\\
\mathbb{H}_n\big(\ell_{\lambda,\alpha}\big)&=&\{(z,t)\in
\mathbb{H}_n; (U_\lambda +\frac{\alpha}{2} z\times z,\alpha
z,\alpha)=(U_\lambda,0,\alpha)\}\\&=&\{0\}\times\mathbb{R}.
\end{eqnarray*}
It follows that
$G\big(\ell_{\lambda,\alpha}\big)=K\big(\ell_{\lambda,\alpha}\big)\ltimes\mathbb{H}_n\big(\ell_{\lambda,\alpha}\big)$.
According to Lipsman [13], the representation
$\pi_{(\lambda,\alpha)}$ is equivalent to the representation of
$G$ obtained by holomorphic induction from the linear form
$\ell_{\lambda,\alpha}$. Now, for an irreducible unitary
representation $\tau_{\mu}$ of $K$ with highest weight
$\mu$, we take the linear functional
$\ell_{\mu}:=(U_\mu,0,0)$ of $ \mathfrak{g}^* $ where
$$U_\mu=
\left( \begin{array}{ccc}
i\mu_1& \ldots & 0\\
\vdots & \ddots & \vdots\\
0 & \ldots & i\mu_n
\end{array}\right).$$
which is clearly aligned and admissible. Hence, the representation of $G$
obtained by holomorphic induction from the linear functional
$\ell_{\mu}$ is equivalent to the representation
$\tau_\mu$.\vspace{0,2cm}
We denote by $ \mathcal{O}^{G}_\mu$ the coadjoint orbit of $\ell_\mu$ and by $\mathcal{O}_{(\lambda,\alpha)}^G$ the coadjoint orbit associated to the linear form $\ell_{\lambda,\alpha}$. Let $\mathfrak{g}^{\ddagger}$ be the set of all admissible linear forms of $G$. The orbit space $\mathfrak{g}^{\ddagger}/G$ is called the space of admissible coadjoint orbits of $G$. The set of all coadjoint orbits $\mathcal{O}_{(\lambda,\alpha)}^G$ turns out to be the subspace of generic admissible coadjoint orbits of $G$.\\

Let $\tau_{\lambda}$ be an irreducible unitary representation of the unitary group $K=U(n)$ with highest weight $\lambda = (\lambda_{1},...,\lambda_{n}) \in \mathbb{Z}^{n}$.
Recall that the irreducible representations of $G = K\ltimes\mathbb{H}_{n}$ that come from an infinite dimensional irreducible representation $\sigma _{\alpha} \in \widehat{\mathbb{H}_{n}}$, $\alpha \in\mathbb{R}^{*} $, are of the form $\pi _{(\lambda,\alpha)}$ with
\begin{equation*}
\pi_{(\lambda,\alpha)}(A,z,t)=\tau_\lambda(A)\otimes\sigma_\alpha(z,t)\circ
W_\alpha(A)
\end{equation*}
for $(A,z,t)\in G$. Here $W_\alpha$ denotes the natural representation of $K$ on the ring $\mathbb{C}[z_1,...,z_n]$ of holomorphic polynomials on $\mathbb{C}^n$, given by $$(A.p)((z_1,...,z_n)^T)=p(A^{-1}(z_1,...,z_n)^T).$$ The space $\mathbb{C}[z_1,...,z_n]$ decomposes under the action of $K$ as $$\mathbb{C}[z_1,...,z_n]=\sum_{k=0}^{\infty}\mathbb{C}_k[z_1,...,z_n]$$ where $\mathbb{C}_k[z_1,...,z_n]$ denotes the space of homogeneous polynomials of degree $k$, thus we have $W_{\alpha}=\displaystyle\bigoplus_{k\in\mathbb{N}}\tau_{\alpha,k}$ where $\tau_{\alpha,k}$ is the representation of $K$ on $\mathbb{C}_k[z_1,...,z_n]$. Consider now an irreducible unitary representation $\tau_\mu$ of $K$ with highest weight $\mu$. The multiplicity of $\tau_\mu$ in the representation $\pi_{(\lambda,\alpha)}$ is given by
\begin{eqnarray*}
m(\pi_{(\lambda,\alpha)},\tau_{\mu})&=&\textit{mult}(\pi_{(\lambda,\alpha)}\big\vert_{K},\tau_\mu)\\
&=&\textit{mult}(\tau_\lambda\otimes W_\alpha,\tau_\mu)\\
&=&\textit{mult}(\bigoplus_{k\in\mathbb{N}}\tau_{\lambda}\otimes\tau_{\alpha,k},\tau_\mu).
\end{eqnarray*}

%%%%%%%%%%%%%%%%%%%%%%%%%%%%%%%%%%%%%%%%%%%%%%%%%%%%%%%%%%%%%%%%%%%%%%%%%%%%%%%%%%%%%%%%%%%%%%%%%%%%%%%%%%%%%%%%%%%%%%%%%%%%
\subsection{Corwin-Greenleaf multiplicity function for $U(n)\ltimes\mathbb{H}_{n}$}
%%%%%%%%%%%%%%%%%%%%%%%%%%%%%%%%%%%%%%%%%%%%%%%%%%%%%%%%%%%%%%%%%%%%%%%%%%%%%%%%%%%%%%%%%%%%%%%%%%%%%%%%%%%%%%%%%%%%%%%%%%%
We continue to use the notation of the previous sections. Fix $\alpha$ a nonzero real. Let $\pi_{(\lambda,\alpha)}\in\widehat{G}$ and $\tau_{\mu}\in\widehat{K}$ be as before. To these unitary representations, we attach respectively the generic coadjoint orbit $\mathcal{O}_{(\lambda,\alpha)}^{G}$ and the coadjoint orbit $\mathcal{O}_{\mu}^{K}$. Here $\mathcal{O}_{\mu}^{K}$ is the orbit in $\mathfrak{k}^{*}$ through $U_{\mu}$ i.e., $\mathcal{O}_{\mu}^{K}=\operatorname{Ad}_{K}^{*}(K)U_{\mu}$.
Now, we turn our attention to the multiplicity $m(\pi_{(\lambda,\alpha)},\tau_{\mu})$ of $\tau_{\mu}$ in the restriction of $\pi_{(\lambda,\alpha)}$ to $K$, we shall prove the following result:

\begin{teo} We have $$m(\pi_{(\lambda,\alpha)},\tau_{\mu})\neq 0\Rightarrow n(\mathcal{O}^{G}_{(\lambda,\alpha)},\mathcal{O}^{K}_{\mu})\neq
0.$$
\end{teo}
\textbf{Proof.}
Denote by $\tau_{\alpha,k}=\tau_{(0,...,0,-k)}$ the irreducible representation of $K$ on $\mathbb{C}_k[z_1,...,z_n]$ with highest weight $(0,...,0,-k)\in\mathbb{Z}^n$. Then, we have
\begin{eqnarray*}
\pi_{(\lambda,\alpha)}\big\vert_{K}&=&\tau_{\lambda}\otimes W_{\alpha}\\
&=&\tau_{\lambda}\otimes \bigoplus_{k\in\mathbb{N}}\tau_{(0,...,0,-k)}\\
&=&\bigoplus_{k\in\mathbb{N}}\tau_{\lambda}\otimes \tau_{(0,...,0,-k)}.
\end{eqnarray*}
Consider again the set $\mathcal{F}_{\mu}=\left\lbrace z\in\mathbb{C}^n;U_{\lambda}+\frac{\alpha}{2}z\times z\in\mathcal{O}_{\mu}^K\right\rbrace $.
Now, assume that $m(\pi_{(\lambda,\alpha)},\tau_{\mu})\neq 0$. Then there exists $k\in\mathbb{N}$ such that $$\tau_{\mu}\subset\tau_{\lambda}\otimes\tau_{(0,...,0,-k)}$$ hence $$\mathcal{O}_{\mu}\subset\mathcal{O}_{\lambda}+\mathcal{O}_{(0,...,0,-k)}$$
So, there exists $C\in U(n)$ such that $$U_{\lambda}+CU_{(0,...,0,-k)}C^{-1}\in\mathcal{O}_{\mu}$$
Let $z=C(0,...,0,r)^t$ with 
$$
r= \left\lbrace\begin{array}{cc}
i\displaystyle\sqrt{\frac{2k}{\alpha}} & \hspace{0.2cm}\mbox{if}\,\, \alpha>0,\\
\\
\displaystyle\sqrt{\frac{-2k}{\alpha}} & \hspace{0.2cm}\mbox{if}\,\,\alpha<0.

\end{array}
\right.$$
Therefore, we have $\frac{\alpha}{2}z\times z=CU_{(0,...,0,-k)}C^{-1}$. It follows that $\mathcal{F}_{\mu}\neq \emptyset$, and then $n(\mathcal{O}^{G}_{(\lambda,\alpha)},\mathcal{O}^{K}_{\mu})\neq 0$.
$\hfill\square$\vspace{0,2cm}\\

The converse of this theorem is false in general if we take for example $\lambda=(-1,...,-1)$ and $\mu=(0,...,0,-1)$ we will see in the last theorem that $n(\mathcal{O}^{G}_{(\lambda,\alpha)},\mathcal{O}^{K}_{\mu})\neq 0$ (see Theorem 4) but \begin{eqnarray*}
\tau_{\lambda}\otimes W_{\alpha}&=&\bigoplus_{k\in\mathbb{N}}\tau_{\lambda}\otimes \tau_{(0,...,0,-k)}\\
&=&\bigoplus_{k\in\mathbb{N}}\tau_{(-1,...,-1,-1-k)}.
\end{eqnarray*}
Therefore $\tau_{\mu}=\tau_{(0,...,0,-1)}\nsubseteq\tau_{\lambda}\otimes W_{\alpha}$ and then $m(\pi_{(\lambda,\alpha)},\tau_{\mu})= 0$.\\

In the remainder of this paper, we give two situation where the Corwin-Greenleaf multiplicity function is less than one and discuss the relationship between $n(\mathcal{O}^{G}_{(\lambda,\alpha)},\mathcal{O}^{K}_{\mu})$ and $m(\pi_{(\lambda,\alpha)},\tau_{\mu})$. For some particular dominant weight $\mu$, we shall prove in the first situation that $m(\pi_{(\lambda,\alpha)},\tau_{\mu})$ coincides with $n(\mathcal{O}^{G}_{(\lambda,\alpha)},\mathcal{O}^{K}_{\mu})$, but in the second situation we have $m(\pi_{(\lambda,\alpha)},\tau_{\mu})\neq n(\mathcal{O}^{G}_{(\lambda,\alpha)},\mathcal{O}^{K}_{\mu})$ .\\

Let us first fix some notation that we will use later. Let $\lambda=(\lambda_1,...,\lambda_n), \mu=(\mu_1,...,\mu_n)\in\mathbb{Z}^n $ such that $\lambda_1\geq...\geq\lambda_n$ and $\mu_1\geq...\geq\mu_n$. To these dominant weights of $K$ we attach the matrix $B_{\lambda,\mu}$ and the vector $V_{\lambda,\mu}$ defined as follows $$B_{\lambda,\mu}=\left( \prod_{k=1,k\neq j}^n(\mu_i-\lambda_k)\right)_{1\leq i,j\leq n}\hspace{0.2cm}\text{and} \hspace{0.2cm} V_{\lambda,\mu}=\left( \prod_{k=1}^n(\mu_1-\lambda_k),...,\prod_{k=1}^n(\mu_n-\lambda_k)\right)^T. $$Now , we are in position to prove
\begin{teo} Let $n\geq 2$. Assume that $\lambda$ is strongly dominant weight of $K$. Then for any dominant weight $\mu$ of $K$ such that $B_{\lambda,\mu}$ is invertible we have $$n(\mathcal{O}^{G}_{(\lambda,\alpha)},\mathcal{O}^{K}_{\mu})\leq 1.$$
\end{teo}
\textbf{Proof.} Let $\lambda=(\lambda_1,...,\lambda_n)$ be a strongly dominant weight of $K$. We shall denote by $H_{\lambda}$ the stabiliser of $U_{\lambda}$ in $K$. Assume that $n(\mathcal{O}^{G}_{(\lambda,\alpha)},\mathcal{O}^{K}_{\mu})\neq 0$ for some dominant weight $\mu$ of $K$. Then there exists $z\in\mathbb{C}^n$ such that $U_{\lambda}+\frac{\alpha}{2}z\times z=AU_{\mu}A^*$ for some $A\in K$. For all $x\in\mathbb{R}$, we have $$det(U_{\lambda}+\frac{\alpha}{2}z\times z -ix\mathbb{I})=(-i)^nP(x)$$ where $P$ is the unitary polynomial of degree $n$ given by $$P(x)=\prod_{i=1}^n(x-\lambda_i)-\frac{\alpha}{2} \sum_{j=1}^n\prod_{i=1,i\neq j}^n(x-\lambda_i)|z_j|^2.$$ Therefore we have $P(\mu_k)=0$ for $k=1,...,n$. It follows that $$V_{\lambda,\mu}=\frac{\alpha}{2}B_{\lambda,\mu}(|z_1|^2,...,|z_n|^2)^T$$ Consider again the set $\mathcal{F}_{\mu}=\left\lbrace z\in\mathbb{C}^n, U_{\lambda}+\frac{\alpha}{2}z\times z\in\mathcal{O}_{\mu}^K\right\rbrace$. Hence $$\mathcal{F}_{\mu}=\left\lbrace z\in\mathbb{C}^n, (|z_1|^2,...,|z_n|^2)^T=\frac{2}{\alpha}B_{\lambda,\mu}^{-1}V_{\lambda,\mu}\right\rbrace .$$ Since
$H_{\lambda}=\mathbb{T}^n$ the $n$-dimensional torus, we conclude that $n(\mathcal{O}^{G}_{(\lambda,\alpha)},\mathcal{O}^{K}_{\mu})= 1$.
$\hfill\square$\vspace{0.2cm}
\begin{corol}
Let $n\geq 2$. Assume that $\lambda=(\lambda_1,...,\lambda_n)$ is strongly dominant weight of $K$ and $\mu=(\lambda_1,...,\lambda_{n-1},\lambda_n-k)$ for some $k\in\mathbb{N}$. Then we have $$m(\pi_{(\lambda,\alpha)},\tau_{\mu})= n(\mathcal{O}^{G}_{(\lambda,\alpha)},\mathcal{O}^{K}_{\mu})$$
\end{corol}
\textbf{Proof.}
Let $\lambda=(\lambda_1,...,\lambda_n)$ be a strongly dominant weight of $K$. Suppose that $\mu=(\lambda_1,...,\lambda_{n-1},\lambda_n-k)$ for some $k\in\mathbb{N}$, then $B_{\lambda,\mu}$ is invertible, therefore $n(\mathcal{O}^{G}_{(\lambda,\alpha)},\mathcal{O}^{K}_{\mu})\leq 1$.
Since $\pi_{(\lambda,\alpha)}\big\vert_{K}=\displaystyle\bigoplus_{k\in\mathbb{N}}\tau_{(\lambda_1,...,\lambda_{n-1},\lambda_n-k)}$ then $m(\pi_{(\lambda,\alpha)},\tau_{\mu})=1$ and by the theorem 2 we deduce that $$m(\pi_{(\lambda,\alpha)},\tau_{\mu})= n(\mathcal{O}^{G}_{(\lambda,\alpha)},\mathcal{O}^{K}_{\mu}).$$
$\hfill\square$\vspace{0,2cm}\\

Concluding this section, let us prove the following result:

\begin{teo} Let $n\geq 2$. If the dominant weight $\lambda=(\lambda_1,...,\lambda_n)$ of $K$ satisfies $\lambda_1=...=\lambda_n=a$ for some $a\in\mathbb{Z}$, then for any dominant weight $\mu$ of $K$ with $\mu\neq\lambda$ we have
$$n(\mathcal{O}^{G}_{(\lambda,\alpha)},\mathcal{O}^{K}_{\mu})\leq 1$$ Moreover, $n(\mathcal{O}^{G}_{(\lambda,\alpha)},\mathcal{O}^{K}_{\mu})\neq0$ if and only if $\mu$ is of the form\vspace{0.2cm}\\
\textbf{Case 1:} if $\alpha>0$ then $\mu=(\underbrace{b,...,b}_p,\underbrace{a,...,a}_q)\in\mathbb{Z}^n$, $p+q =n$, $b\in \mathbb{Z}$ with $b>a.$\\
\textbf{Case 2:} if $\alpha<0$ then $\mu=(\underbrace{a,...,a}_p,\underbrace{b,...,b}_q)\in\mathbb{Z}^n$, $p+q =n$, $b\in \mathbb{Z}$ with $a>b.$\\
Consequently, if $\mu_{n-1}\neq a$ and $n(\mathcal{O}^{G}_{(\lambda,\alpha)},\mathcal{O}^{K}_{\mu})\neq 0$ then $m(\pi_{(\lambda,\alpha)},\tau_{\mu})\neq n(\mathcal{O}^{G}_{(\lambda,\alpha)},\mathcal{O}^{K}_{\mu})$.
\end{teo}
\textbf{Proof.}
Let $\lambda=(\lambda_1,...,\lambda_n)$ be a dominant weight of $K$ such that $\lambda_1=...=\lambda_n=a$ with $a\in\mathbb{Z}$. Assume that $n(\mathcal{O}^{G}_{(\lambda,\alpha)},\mathcal{O}^{K}_{\mu})\neq 0$ for some dominant weight $\mu$ of $K$. Then there exists $z\in\mathbb{C}^n$ such that $U_{\lambda}+\frac{\alpha}{2}z\times z=AU_{\mu}A^*$ for some $A\in K$. For all $x\in\mathbb{R}$, we have $$det(U_{\lambda}+\frac{\alpha}{2}z\times z -ix\mathbb{I})=(-i)^nP(x)$$ with $$P(x)=(x-a)^{n-1}\left( x-a- \frac{\alpha}{2} \sum_{j=1}^n|z_j|^2\right). $$ Then we have $P(\mu_k)=0$ for $k=1,...,n$. It follows that $$\left\{\begin{array}{lll}
\mu_k = a \\
\mbox{or}\\
\mu_k \neq a \hspace{0.2cm}\mbox{and}\hspace{0.2cm}\displaystyle \mu_k = a+ \frac{\alpha}{2}\sum_{j=1}^n|z_j|^2.
\end{array}
\right.$$ Since $\mu\neq\lambda$ then there exists $1\leq k\leq n$ such that $\mu_k\neq a.$\\
\textbf{\underline{Case $\alpha>0$}} : Let $p=max\{1\leq k\leq n,\mu_k\neq a\}$
 then $$\mu_p=a+ \frac{\alpha}{2}\sum_{j=1}^n|z_j|^2>a.$$ Since $\mu_1\geq...\geq\mu_p\geq...\geq\mu_n$, we obtain $$\mu=(\underbrace{b,...,b}_p\underbrace{a,...,a}_q) \hspace{0.2cm}\text{with}\hspace{0.2cm} b=a+ \frac{\alpha}{2}\sum_{j=1}^n|z_j|^2.$$
Consider again the set $\mathcal{F}_{\mu}=\left\lbrace z\in\mathbb{C}^n, U_{\lambda}+\frac{\alpha}{2}z\times z\in\mathcal{O}_{\mu}^K\right\rbrace $ then $$\mathcal{F}_{\mu}=\left\lbrace z\in\mathbb{C}^n, \sum_{j=1}^n|z_j|^2=(b-a) \frac{2}{\alpha}\right\rbrace .$$ Since $H_{\lambda}=K$ we can deduce that $n(\mathcal{O}^{G}_{(\lambda,\alpha)},\mathcal{O}^{K}_{\mu})=1.$\\
\textbf{\underline{Case $\alpha<0$}} : Let $l=min\{1\leq k\leq n,\mu_k\neq a\}$ then $$\mu_l=a+ \frac{\alpha}{2}\sum_{j=1}^n|z_j|^2<a.$$ Hence $$\mu=(\underbrace{a,...,a}_p\underbrace{b,...,b}_q)\hspace{0.2cm} \text{with}\hspace{0.2cm} b=a+ \frac{\alpha}{2}\sum_{j=1}^n|z_j|^2, p=l-1$$ and so $n(\mathcal{O}^{G}_{(\lambda,\alpha)},\mathcal{O}^{K}_{\mu})=1.$\\
Now, Suppose that $\mu_{n-1}\neq a$, if $\alpha>0$ we get $\mu=(b,...,b,a)\in\mathbb{Z}^{n}$ with $b>a$ and if $\alpha<0$, $\mu=(\underbrace{a,...,a}_p,\underbrace{b,...,b}_q)\in\mathbb{Z}^n$ with $a>b$ and $q\geq 2$. Since $\pi_{(\lambda,\alpha)}\big\vert_{K}=\displaystyle\bigoplus_{k\in\mathbb{N}}\tau_{(a,...,a,a-k)}$ then $m(\pi_{(\lambda,\alpha)},\tau_{\mu})=0$ and hence $m(\pi_{(\lambda,\alpha)},\tau_{\mu})\neq n(\mathcal{O}^{G}_{(\lambda,\alpha)},\mathcal{O}^{K}_{\mu})$.\vspace{0.2cm}\\
This completes the proof of the theorem.
$\hfill\square$\vspace{0,2cm}\\

%%%%%%%%%%%%%%%%%%%%%%%%%%%%%%%%%%%%%%%%%%%%%%%%%%%%%%%%%%%%%%%%%%%%%%%%%%%%%%%%%%%%%%%%%%%%%%%%%%%%%%%%%%%%%%%%%%%%%%%%%%%%%%%%%%%%%%%%%%%%%%%%%%%%%%%%%%%%%%%%%
\section*{Acknowledgment}
The authors would like to thank the referee for his / her careful reading of our paper and for remarks improving the article.

\small{

\vspace{0,5cm}

Department of Mathematics, Faculty of Sciences at Sfax, University of Sfax, Route de Soukra, B. P. 1171, 3000-Sfax, Tunisia\\
E-mail adress: majdi.benhalima@yahoo.fr\\

Department of Mathematics, Preparatory Institute for Engineering Studies at Gafsa, University of Gafsa, El Khayzorane street-Zaroug, 2112-Gafsa, Tunisia\\
E-mail adresses: anis.messaoud@ipeig.rnu.tn

\end{document}